\documentclass[a4paper]{article}

\usepackage{epsf}
\usepackage{color}
\usepackage{epsfig}
\usepackage{graphicx}
\usepackage{amsmath}
\usepackage{latexsym}
\usepackage{amssymb}
\usepackage{a4}

\newcommand{\R}{\mathbb{R}}

\newcommand{\e}{{\rm e}}

\newcommand{\C}{\mathbb{C}}

\providecommand{\keywords}[1]
{
  \small	
  \textbf{\textit{Keywords---}} #1
}

\begin{document}

\title{Computing the matrix sine and cosine simultaneously with a reduced number of products}

\author{Muaz  Seydao\u{g}lu$^{1,3}$ \footnote{ E-mail: m.seydaoglu@alparslan.edu.tr}, Philipp Bader$^{2}$\footnote{  bader@uji.es}, Sergio Blanes$^{3}$ \footnote{serblaza@imm.upv.es} Fernando Casas$^{4}$\footnote{casas@uji.es}  \\
        \small $^{1}$Faculty of Art and Science, Department of Mathematics
49100 Mus, Turkey. \\
        \small $^{2}$Departament de Matem\`atiques, Universitat Jaume I, 12071 Castell\'on, Spain.\\
				\small $^{3}$Instituto de Matem\'atica Multidisciplinar,
  Universitat Polit\`ecnica de Val\`{e}ncia, E-46022  Valencia.\\
				\small $^{4}$IMAC and Departament de Matem\`atiques, Universitat Jaume I, 12071 Castell\'on.\\
}

\date{}
\maketitle

%
%
%

\begin{abstract}
A new procedure is presented for computing the matrix cosine and sine simultaneously by means of Taylor polynomial approximations. These
are factorized so as to reduce the number of matrix products involved. Two versions are developed to be used in single and double precision
arithmetic. The resulting algorithms are more efficient than schemes based on Pad\'e approximations for a wide range of norm matrices.
\end{abstract}\hspace{10pt}


\section{Introduction} \label{Introduction}

Many dynamical systems are modeled by differential equations in which finding closed solutions is not possible and so one has to compute approximating solutions. These differential equations usually preserve some underlying geometric structure
which reflects the qualitative nature of the phenomena they describe. It is then relevant that the approximations share with the
exact solution of the differential equation these qualitative properties to render a description. 
The
design and analysis of numerical integrators preserving some of these geometric structures constitutes the realm of Geometric Numerical Integration (GNI), an active and interdisciplinary research area and the subject of intensive development during the last decades \cite{blanes16aci,hairer06gni,iserles00lgm1,leimkuhler04shd,lubich08fqt,sanzserna94nhp}.

Exponential integrators can be considered as a class of GNIs tailored to stiff and oscillatory equations \cite{blanes15aea,blanes09tme,hochbruck97oks,hochbruck10ein,iserles00lgm}. For large systems of equations these schemes usually require to compute the action of the exponential of a matrix on a vector \cite{hochbruck97oks,hochbruck10ein}. However, for problems of moderate size it may be more appropriate to compute directly the exponential of the matrices involved. 

When the problem is oscillatory, very often the formal solution involves both the sine and cosine of a matrix. Thus, for example, consider
the Schr\"odinger equation in quantum mechanics,
\[
  i \frac{d \psi}{dt} = {\cal H}(t) \psi, \qquad \psi(t_0)=\psi_0,
\]
where ${\cal H}(t)$ is a Hermitian operator and $\psi$ is a complex wave function. A usual procedure to get numerical approximations involves first a 
spatial discretisation or working on a finite dimensional representation. In any event, one ends up with a matrix equation with a similar structure,
\[
  i  \frac{du}{dt} = A u, \qquad u(t_0)=u_0\in \C^N.
\]
If $A$ is a real and constant matrix, the unitary evolution operator is given by
\begin{equation} \label{eq:unitary}
   U(t)= \e^{-itA} = \cos(tA) -i\sin(tA).	
\end{equation}

There are different techniques to compute efficiently the exponential of a matrix \cite{preprint,bader19ctm,bader15tss,higham10cma,moler03ndw,sastre19ftp,sastre19btc,sidje98eas}. However, using any of these general algorithms to approximate the unitary matrix $\e^{-itA}$ in \eqref{eq:unitary} involves products of complex matrices 
making them computationally expensive. Alternatively, we propose an efficient procedure to compute the matrix sine and cosine that only involves a small number of products of real matrices. The algorithm is used in combination with the squaring as
\[
\cos(2A)=2\cos^2(A)-I=I-2\sin^2(A), \qquad \sin(2A)=2\sin(A)\cos(A).
\]
In this way, it only requires two products per squaring (instead of four products when considering the square of complex matrices), thus
making the overall procedure more efficient.

There are other examples where the computation of the sine and cosine of a matrix can be of interest. For example, for wave equations given by the generic second order system
\[
  y'' + A y = f(y,t),
\]
with $y\in \R^N$, exponential integrators frequently require to solve separately the linear homogeneous problem
\begin{equation}  \label{eq.2a}
  y'' + A y = 0, \qquad y(0)=y_0, \quad y'(0)=y'_0.
\end{equation}
Writing (\ref{eq.2a}) as a first order system, the solution is given by
\begin{equation}
	\left( \begin{array}{c} y(t) \\ y'(t)
	\end{array}	\right) =
	\e^{tM}
	\left( \begin{array}{c} y_0 \\ y'_0
	\end{array}	\right), 
	\qquad \mbox{with} \quad
	M= \left( \begin{array}{cc} 0 & I \\ -A & 0
	\end{array}	\right)
\end{equation}
and
\begin{equation}\label{eq:symplectic}
\begin{array}{rcl}
	\e^{tM} &=&
	\left( \begin{array}{cc} 
	\cos(t\sqrt{A}) & (\sqrt{A})^{-1}\sin(t\sqrt{A}) \\ -\sqrt{A}\sin(t\sqrt{A}) & \cos(t\sqrt{A})
	\end{array}	\right) 
	 \equiv \left( \begin{array}{cc} 
	c(t^2{A}) & s(t,{A}) \\ -{A}s(t,{A}) & c(t^2{A})
	\end{array}	\right).
	\end{array}
\end{equation}

 Notice that the dimension of 
$M$ is twice the dimension of $A$ and so the cost of matrix-matrix multiplications grows, in general, by a factor of eight.

On the other hand, a closer look to the functions to be approximated clearly indicates that the same algorithm used to evaluate the matrix sine and cosine for the unitary matrix (\ref{eq:unitary}) should not be used directly since it requires computing first the square root of the matrix, $B=\sqrt{A}$, in addition to a multiplication and an inversion of this matrix. As a matter of fact, an efficient approximation to the exponential \eqref{eq:symplectic} was already presented in \cite{bader17sif}. We propose in this case an improved algorithm based on a modification of the methods to compute the matrix sine and cosine with the goal of computing simultaneously the functions $c(t^2A) \equiv \cos(\sqrt{t^2A})$ and 
$s(t,A) \equiv (\sqrt{A})^{-1}\sin(\sqrt{t^2A})$. For the double angle we will take into account that
\[
c(4t^2A)=2c^2(t^2A)-I, \qquad s(2t,A)=2s(t,A)c(t^2A),
\]
thus requiring only two products per squaring. Notice that we do not use the property $\cos(2A)=I-2\sin^2(A)$ since the function $\sin(A)$ is not computed in this case.

In summary, the purpose of this paper consists in developing algorithms that allow one to compute 
$\cos(A)$ and $\sin(A)$ or $c(t^2A)$ and $s(t,A)$ simultaneously and providing full accuracy up to single or double precision with a reduced 
computational cost. Thus, in particular, 
we propose an algorithm that, with only four products, approximates $\cos(A)$ with an error of order ${\cal O}(A^{17})$, and with two extra products it also approximates $\sin(A)$ with an error of order ${\cal O}(A^{18})$. The same procedure allows one,  
with one extra product (seven products in total), to approximate $\cos(A)$ and $\sin(A)$ with errors of order ${\cal O}(A^{25})$ and ${\cal O}(A^{24})$, respectively. 

Although one can find in the literature several algorithms to compute $\cos(A)$ (see \cite{sastre19ftp} and references therein), only few of them
are designed to do so in a simultaneous way (see \cite{almohy14naf} and references therein). As our analysis shows and several numerical
examples confirm, the technique we propose here outperform all of them.

\section{The algorithms}

The search of fast algorithms for  evaluating matrix polynomials has received considerable interest in the recent literature \cite{preprint,bader19ctm,lei92afa,liang04fmf,paterson73otn,sastre18eeo,vanloan79ano,westreich89etm}. We next briefly summarize 
 how to approximate the matrix sine and cosine functions by means of certain polynomials involving a reduced number of matrix products.
 This reduction essentially follows the same approach used in \cite{blanes02hoo} to minimise the number of commutators 
 appearing in different Lie-group integrators and was successfully adapted to the Taylor expansion of the exponential 
 matrix in \cite{preprint} and especially in \cite{bader19ctm}. 

Generally speaking, the strategy consists first  in elaborating a recursive procedure to compute the polynomial approximating the matrix
cosine with the minimum number of products and then these same products are used to approximate the matrix sine as accurately as possible
in the cheapest possible way.

Clearly, the most economic way to construct polynomials of degree $2^k$ is by applying the following sequence, which requires the
evaluation of only $k$ products. First we form the intermediate matrices
\begin{align}
A_0 &= I, \quad A_1 = A, \nonumber  \\
A_2 &= z_{2,0}I+z_{2,1}A_1+(x_1I+x_2A_1)(x_3 I+x_4 A_1), \nonumber \\
A_4 &= z_{4,0}I+z_{4,1}A_1+z_{4,2}A_2+(x_5I+x_6 A_1+x_7A_2)(x_8I+x_9A_1+x_{10}A_2),  \nonumber
\\
A_8 &= \sum_{k=0}^3z_{8,2^{k-1}}A_{2^{k-1}}+ (x_{11}I+\cdots+x_{14}A_4)(x_{15}I+\cdots+x_{18}A_4), \label{Gen_algorithm} \\
 & \vdots  \nonumber
\end{align}
and finally we take
\[
P_{2^k} = A_{2^k}.
\]
Here the indices in $A$, $A_{2^k}$, are chosen to indicate the highest attainable power, i.e., $A_{2^k}=\mathcal{O}(A^{2^k})$. Of course, 
there are many redundancies in the coefficients since some of them can be absorbed by others.

It is a simple exercise to check that any polynomial of degree up to four can be computed with two products, whereas polynomials up to degree eight  can be computed with only three products. This does not mean, however, that all such polynomials can be written with just three products. 
This is the case, in particular, of $P_7(A)=A^7$, as can be readily seen. When a given polynomial cannot be reproduced by following the previous
approach, new terms can be incorporated. Thus, in particular 
\begin{align}
A_0 &= I, \ \ A_1 = A, \ \ A_2=A^2, \ \ A_3=A A_2\nonumber  \\
A_6 &= B_{3,1}+B_{3,2}B_{3,3}, \qquad 
  B_{3,i}=\sum_{k=0}^3 x_{i,k}A_k \label{Mod_algorithm}
\\
 & \vdots  \nonumber
\end{align}
and this generalises the procedure.

We use this technique in the sequel to approximate first $\cos(A)$ and $\sin(A)$ simultaneously with the minimum number of products, 
and then we apply the same procedure to $c(t^2A)$ and $s(t,A)$.

\subsection{Computing $\cos(A)$ and $\sin(A)$ simultaneously}

Let us denote by
\[
  T_{2m}^c=\sum_{k=0}^m\frac{(-1)^k(A^2)^{k}}{(2k)!}, \qquad
  T_{2m+1}^s=A \, \sum_{k=0}^m\frac{(-1)^k(A^2)^{k}}{(2k+1)!}
\]
the Taylor polynomial approximations of $\cos(A)$ and $\sin(A)$ up to order $2m$ and $2m+1$ in $A$, respectively, and by
$T_{2m+1,\ell}^s$ with $\ell>2m+1$, any polynomial of degree $\ell$ such that 
\[
   T_{2m+1,\ell}^s=T_{2m+1}^s+{\cal O}(A^{2m+2}).
\]   

\paragraph{$k=3$ products.}
This constitutes a trivial problem, but it nevertheless illustrates the general procedure. With two products we can compute  $T_4^c$:
\begin{equation}  \label{Algorithm3cA}
\begin{aligned}
   A_2 & = A^2,\\
   A_4 & = A^4,\\
   T_4^c(A)  & = I - \frac1{2!} A_2+ \frac1{4!}A_4,
\end{aligned}  
\end{equation}
and with one extra product we can get
\begin{equation}  \label{Algorithm3sA}
   T_5^s(A) = A(I - \frac1{3!} A_2+ \frac1{5!}A_4). 
\end{equation}

\paragraph{$k=4$ products.}
 With three products we can compute  $T_8^c$:
\begin{equation}  \label{Algorithm4cA}
\begin{aligned}
   A_2 & = A^2,\\
   A_4 & = A_2^2,\\
   A_8 & = A_4\left(-\frac1{6!} A_2+\frac1{8!} A_4\right),\\
   T_8^c(A)  & =  I - \frac1{2!} A_2 + \frac1{4!} A_4+ A_8. 
\end{aligned}  
\end{equation}
With one extra product we can approximate the matrix sine, but only up to order seven as follows
\begin{equation}  \label{Algorithm4sA}
   T_{7,9}^s(A)  = A\left(I - \frac1{3!} A_2+ \frac1{5!}A_4 + \frac{6!}{7!}A_8\right). 
\end{equation}
According with the previous notation, $T_{7,9}^s(A)  = T_{7}^s(A) + \mathcal{O}(A^8)$. 

The order of approximation of the matrix sine can be increase up to order nine by incorporating one extra product as follows:
\begin{equation}  \label{Algorithm5sA}
\begin{aligned}
   A_8 & = A_4\left(-\frac1{7!} A_2+\frac1{9!} A_4\right),\\
   T_9^s(A)  & =  A\left(I - \frac1{3!} A_2 + \frac1{5!} A_4+ A_8\right). 
\end{aligned}  
\end{equation}

\paragraph{$k=6$ products.}
 The following scheme allows one to express $T_{16}^c(A)$ with only four products:
\begin{equation}  \label{Algorithm6cA}
\begin{aligned}
   A_2 & = A^2,\\
   A_4 & = A_2^2,\\
   A_8 & = A_4(x_1 A_2+x_2 A_4),\\
   A_{16} & = (x_3 A_4+A_8)(x_4I+x_5A_2+x_6 A_4+x_7A_8),  \\
   T_{16}^c(A)  & = I -\frac12 A_2 +x_{8} A_4+ A_{16}. 
\end{aligned}  
\end{equation}
In fact, we get two families of solutions depending on a free parameter, $x_1$, which is chosen to (approximately) minimize the 1-norm of the vector of parameters $(x_1,\ldots, x_8)$. This results in

\begin{equation} \label{eq:cosm8}
\begin{array}{lll}
x_1=\displaystyle \frac{ 7}{500},& 
x_2= \displaystyle -\frac{7 }{60000},&
x_3= \displaystyle \frac{1 }{2500} (-1533  + 7 \sqrt{36681}), \\
x_4= \displaystyle -\frac{5(124581 + 391 \sqrt{36681})}{10594584}, &
x_5= \displaystyle \frac{9775}{10594584},&
x_6=  \displaystyle -\frac{5(1001 + \sqrt{36681})}{508540032},\\
x_7=  \displaystyle \frac{3125}{889945056},& &
x_8 = \displaystyle \frac{1549211 + 3246 \sqrt{36681}}{63063000}.
\end{array}
\end{equation}
Some of  the coefficients are irrational numbers because they correspond to solutions of a nonlinear system of equations.

With two extra products we can approximate the matrix sine up to order $\mathcal{O}(A^{18})$ as follows:
\begin{equation}  \label{Algorithm6sA}
\begin{aligned}
   C_{24} & = (z_5 I + z_5 A_2 +z_{6} A_4+ z_{7} A_8+ z_{8} T_{16}^c(A))A_8,  \\
   T_{17,25}^s(A)  & = A\left(z_0 I + z_1 A_2 +z_{2} A_4+ z_{3} A_8+ z_{4} T_{16}^c(A) + C_{24}\right)
\end{aligned}  
\end{equation}
with
\begin{equation}
\begin{array}{llll}
z_0= \displaystyle \frac{8887}{4794},& 
z_1=-\displaystyle \frac{1897}{3196} ,&
z_2=\displaystyle  \frac{25259}{575280}, \\ 
z_3= -\displaystyle \frac{965093875}{9674368704},&
z_4= -\displaystyle \frac{4093}{4794},&
z_5= \displaystyle \frac{25698275}{29023106112},&\\
z_6= -\displaystyle \frac{3907675}{348277273344},& 
z_7= \displaystyle \frac{11865625}{3656911370112},& 
z_8 = \displaystyle \frac{25}{308756448},
\end{array}
\end{equation}

i.e. it approximates the matrix sine up to a higher order than the matrix cosine.

\paragraph{$k=7$ products.}
 With five products we can compute  $T_{24}^c$:
\begin{equation}  \label{Algorithm7cA}
	\begin{aligned}
    A_2 & = A^2,\\
		A_4&=A_2^2,\\
		A_6&=A_4A_2,\\
		C_1 &= a_{0,1}I+a_{1,1}A_2+a_{2,1}A_4+a_{3,1}A_6,\\
		C_2 &= a_{0,2}I+a_{1,2}A_2+a_{2,2}A_4+a_{3,2}A_6,\\
		C_3 &= a_{0,3}I+a_{1,3}A_2+a_{2,3}A_4+a_{3,3}A_6,\\
		C_4 &= a_{0,4}I+a_{1,4}A_2+a_{2,4}A_4+a_{3,4}A_6,\\
		A_{12} &= C_3 + C_4^2 \\
		A_{24} &= (C_2 + A_{12})A_{12} \\
		T_{24}^c(A) & = C_1 + A_{24}.
	\end{aligned}
	\end{equation}
The best solution we have obtained is:
\begin{equation}  \label{Coefs7cA}
		\begin{array}{ll}
	a_{0,1} =0, &
	a_{1,1} =0, \\
	a_{2,1} = 0.02264979811206039519,&
	a_{3,1} =-0.00013110924142135755, \\
	a_{0,2} = 0.55751443809990408029,&
	a_{1,2} =-0.61577924683458386455,\\
	a_{2,2} = 0.00747198841446687051,&
	a_{3,2} =-0.00003362444420476012,\\
	a_{0,3} = 0.75936877868464999248,&
	a_{1,3} =-0.01560333979813817129,\\
	a_{2,3} = 0.00010936989591908396,&
	a_{3,3} =-1.03893360877457159499\cdot 10^{-6},\\
	a_{0,4} =0,&
	a_{1,4} =-0.039649968743474473091,\\
	a_{2,4} = 0.000155490073503821463,&
	a_{3,4} =-1.126739663071170022488\cdot 10^{-6}.
	\end{array}
	\end{equation}
Although we report here 20 digits for the coefficients, they can be in fact determined with arbitrary accuracy.

With two extra products we can approximate the matrix sine up to order $\mathcal{O}(A^{23})$ as follows:
\begin{equation}  \label{Algorithm7sA}
\begin{aligned}
   C_{48} & = (z_6 I + z_7 A_2 +z_{8} A_4+ z_{9} A_6+ z_{10} A_{12}+ z_{11} T_{24}^c(A))T_{24}^c(A),  \\
   T_{23,49}^s(B)  & = A\left(z_0 I + z_1 A_2 +z_{2} A_4+ z_{3} A_6+ z_{4} A_{12} + z_{5} T_{24}^c(A) + C_{48}\right),
\end{aligned}  
\end{equation}
with 
\begin{equation}  \label{Coefs7sA}
\begin{array}{ll}
z_0=    0.10090808375109885598,&
z_1=   -0.07668753546445299316,\\
z_2=    0.00084924846993243257, &
z_3=   -0.00001220406904464391,\\
z_4=    0.98499703159318860027,&
z_5=   -0.84925233648155398756,\\
z_6=    1,&
z_7=    0.00095544138280925799,\\
z_8=    4.56337109377154270633\cdot 10^{-6},&
z_9=    2.73461259403000427141\cdot 10^{-8},\\
z_{10}= 0.00048550288474842477 &
z_{11}=-4.15891109384923342531\cdot 10^{-7} .
\end{array}
\end{equation}

\subsection{Computing $c(t^2A)$ and $s(t,A)$ simultaneously}

Let us denote by
\[
  P_{m}^c(t^2A)=\sum_{k=0}^m\frac{(-1)^k(t^2A)^{k}}{(2k)!}, \qquad\quad
  P_{m}^s(t,A)=t \, \sum_{k=0}^m\frac{(-1)^k(t^2A)^{k}}{(2k+1)!}
\]
the Taylor expansions of the functions 
\[
  c(t^2{A}) = \cos(\sqrt{t^2A}), \qquad \mbox{ and } \qquad  s(t,{A}) = (\sqrt{A})^{-1}\sin(\sqrt{t^2A})
\]  
up to order $m$ in $A$, respectively,
with $A$ a real matrix. Notice that they are approximations up to order  $2m$ and $2m+1$ in $t$ to the respective functions.
Analogously, we will denote by $P_{m,\ell}^s$, $\ell>m$, any polynomial of degree $\ell$ such that $P_{m,\ell}^s=P_{m}^s+{\cal O}(A^{m+1}).$

Next we show how the previous algorithms to approximate the sine and cosine functions can be adjusted to approximate  $c(t^2A)$ and $s(t,A)$. 
As before, we proceed according with the number of products involved. 

\paragraph{$k=3$ products.}
With two products we can compute  $P_4^c(t^2 A)$:
\begin{equation}  \label{Algorithm83eA}
\begin{aligned}
   B & = t^2A,\\
   B_2 & = B^2,\\
   B_4 & = B^2( -\frac1{6!}B+\frac1{8!}B_2),\\
   P_4^c(B)  & = I - \frac1{2!} B+ \frac1{4!}B_2 + B_4.
\end{aligned}  
\end{equation}
With the same number of products we can also evaluate $P_{3,4}^s(t,A)$,
\begin{equation}  \label{Algorithm83eA}
   P_{3,4}^s(t,A) = t\left(I - \frac1{3!} B+ \frac1{5!}B_2 -\frac{6!}{7!}B_4\right), 
\end{equation}
whereas with one extra product we get
\begin{equation}  \label{Algorithm83eA}
   P_4^s(t,A) = t\left(I - \frac1{3!} B+ \frac1{5!}B_2 + 
	B_2\left(-\frac1{7!}B+\frac1{9!}B_2\right)\right). 
\end{equation}

\paragraph{$k=4$ products.}
 With three products we can compute  $P_8^c(t^2 A)$:
\begin{equation}  \label{Algorithm83eA}
\begin{aligned}
   B & = t^2A,\\
   B_2 & = B^2,\\
   B_4 & = B_2(x_1 B+x_2 B_2),\\
   B_8 & = (x_3 B_2+B_4)(x_4I+x_5B+x_6 B_2+x_7B_4),  \\
   P_8^c(B)  & =  y_0 I + y_1 B +y_{2} B_2+ B_8,
\end{aligned}  
\end{equation}
whose coefficients are the same as those given in \eqref{eq:cosm8}.

With one extra product we can approximate the matrix sine up to order eight as 
\begin{equation}  \label{Algorithm83eA}
\begin{array}{rl}
   C_{12} & = (z_5 I + z_5 B +z_{6} B_2+ z_{7} B_4+ z_{8} P_8^c(B))B_4,  \\
   P_{8,12}^s(t,A) & = t\left(z_0 I + z_1 B +z_{2} B_2+ z_{3} B_4+ z_{4} P_8^c(B) + 
					C_{12}\right),
\end{array}
\end{equation}
with the same values for the coefficients $z_i$ as before.

\paragraph{$k=5$ products.}
 With four products we can compute  $P_{12}^c(t^2 A)$:
\begin{equation}  \label{poly12}
	\begin{aligned}
    B & = t^2A,\\
		B_2&=B^2,\\
		B_3&=B_2B,\\
		D_1 &= a_{0,1}I+a_{1,1}B+a_{2,1}B_2+a_{3,1}B_3,\\
		D_2 &= a_{0,2}I+a_{1,2}B+a_{2,2}B_2+a_{3,2}B_3,\\
		D_3 &= a_{0,3}I+a_{1,3}B+a_{2,3}B_2+a_{3,3}B_3,\\
		D_4 &= a_{0,4}I+a_{1,4}B+a_{2,4}B_2+a_{3,4}B_3,\\
		B_6 &= D_3 + D_4^2 \\
		P_{12}^c(B) & = D_1 + (D_2 + B_6)B_6,
	\end{aligned}
	\end{equation}
with solution for the coefficients $a_{i,j}$ given in \eqref{Coefs7cA}, whereas with
one extra product we can approximate $P^s_{11}(t,A)$ as 
\begin{equation}  \label{Algorithm83eA}
\begin{aligned}
   C_{24} & = (z_6 I + z_7 B +z_{8} B_2+ z_{9} B_3+ z_{10} B_6+ z_{11} P_{12}^c(A))P_{12}^c(B),  \\
   P^s_{11,24}(t,A)  & = t\left(z_0 I + z_1 B +z_{2} B_2+ z_{3} B_3+ z_{4} B_6 + z_{5} P_{12}^c(B) + C_{24}\right),
\end{aligned}  
\end{equation}
with the same coefficients as in \eqref{Coefs7sA}.

\subsection{Pad\'e approximations}

At this point it is useful to briefly review
the schemes presented in \cite{almohy14naf} to compute the matrix sine and cosine simultaneously, since they will be compared in section
\ref{num.exp}  with our own procedure.

The methods presented in \cite{almohy14naf} are based on the identities
\[
  \cos(A)=\frac{\e^{iA}+\e^{-iA}}{2}, \qquad\quad 
  \sin(A)=\frac{\e^{iA}-\e^{-iA}}{2i}, 
\]
and the use of Pad\'e approximations of the exponentials $\e^{i A}$. 
For instance, taking a diagonal Pad\'e of order eight for approximating $\e^{i A}$, i.e. $r_4(iA)=[p_4(-iA)]^{-1} p_4(iA)=\e^{iA}+{\cal O}(A^9)$ 
one gets
\begin{eqnarray}
	s_4 & = & \frac{A\left(I-\frac{11}{8}A^2+\frac{37}{1176}A^4-\frac{1}{70560}A^6\right)}{I+\frac{1}{28}A^2+\frac{3}{3920}A^4+\frac{1}{8}A^6+\frac{1}{2822400}A^8}, \\
	c_4 & = & \frac{I-\frac{13}{28}A^2+\frac{289}{11760}A^4-\frac{19}{70560}A^6+\frac{19}{2822400}A^8}{I+\frac{1}{28}A^2+\frac{3}{3920}A^4+\frac{1}{8}A^6+\frac{1}{2822400}A^8}, 
\end{eqnarray}
where 
\[
  s_4=\sin(A)+{\cal O}(A^9), \qquad
  c_4=\cos(A)+{\cal O}(A^{10}).
\]
It is clear that $s_4,c_4$ can be computed simultaneously with 5 products ($A^2$, $A^4$, $A^6$, $A^8$, and the extra product for the numerator in $s_4$) and the computation of two inverse matrices. Since both denominators are the same, only one $LU$ factorization is necessary. The totals cost is $(7+\frac13)$ products. Notice that the same order (with very similar accuracy as we will see) is obtained with our novel approach at the cost of only 4 products (and a smaller number of matrices need to be stored).

\section{Error analysis}
Next we analyse how to bound the truncation errors of the previously considered Taylor polynomial approximations of order $2m$ and $2\tilde{m}+1$ for cosine and sine functions, respectively. They have the form
\begin{equation}
\begin{aligned} \label{newtyerrsertay}
  & \cos(A)-T_{2m}^{c} = \sum_{k=m+1}^\infty \alpha_{2k}A^{2k}, \qquad 2m\in \left\{4,8,16,24\right\} \\
  & \sin(A)-T_{2\tilde{m}+1,\ell}^{s} = \sum_{k=\tilde{m}+1}^\infty \tilde{\alpha}_{2k+1}A^{2k+1}, \qquad 2\tilde{m}+1\in \left\{5,7,17,23\right\}.
	\end{aligned}
\end{equation}
On the other hand, the truncation errors of the approximations of the cosine and sine functions obtained by using Pad\'{e} approximants for
$\e^{i A}$ \cite{almohy14naf} can be written as 
\begin{equation}  \label{newtyerrserpad}
 \cos(A)-c_{m} =  \sum_{k=m+1}^\infty \gamma_{2k}A^{2k},
  \qquad\quad 
 \sin(A)-s_{m} =  \sum_{k=m}^\infty \hat{\gamma}_{2k+1}A^{2k+1}. 
\end{equation}
Clearly, the series (\ref{newtyerrsertay}) and (\ref{newtyerrserpad}) can be bounded in terms of $\left\|A\right\|$ as
\begin{equation}
\begin{aligned} \label{newtyerrnorm}
 \left\|\cos(A)-T_{2m}^{c}\right\|&\leq& \sum_{k=m+1}^\infty \left|\alpha_{2k}\right|\theta^{2k},
  \qquad\quad 
 \left\|\sin(A)-T_{2\tilde{m}+1,\ell}^{s}\right\|&\leq& \sum_{k=\tilde{m}+1}^\infty \left|\tilde{\alpha}_{2k+1}\right|\theta^{2k+1}, 
\end{aligned}
\end{equation}
and
\begin{equation}
\begin{aligned} \label{paderrnorm}
 \left\|\cos(A)-c_{m}\right\|&\leq& \sum_{k=m+1}^\infty \left|\gamma_{2k}\right|\theta^{2k},
 \qquad\quad 
\left\|\sin(A)-s_{m}\right\|&\leq& \sum_{k=m}^\infty \left|\hat{\gamma}_{2k+1}\right|\theta^{2k+1},
\end{aligned}
\end{equation}
where 
\[
    \theta=\theta(A)= \left\|A\right\|.
\]

We denote by $\theta_{2m}^M$ the largest value of $\theta$ such that the bounds (\ref{newtyerrnorm}), (\ref{paderrnorm})  do not exceed a prescribed accuracy, $u$, for each method $M \equiv T_{2m}^{c}, T_{2\tilde{m}+1}^{s}, c_m, s_m$.  To achieve maximum accuracy, we bound the 
previous forward absolute errors with the unit round off $u = 2^{-53}$, $u = 2^{-24}$ in double and single precision floating-point arithmetic, respectively. We have truncated the series of the corresponding functions
after 150 terms to find $\theta_{2m}^M$.  
The corresponding values for the new Taylor approximations  of the cosine and sine functions  
are collected in Tables \ref{doubleprctable} and \ref{singleprctable}. For completeness, we also include the values of $\theta_{2m}^M$ for
the Pad\'e approximations, as given in \cite{almohy14naf}, and the total number of matrix products corresponding to each procedure $\Pi_{2m}$.
In the case of Pad\'e approximants, we have added the cost of evaluating two inverse matrices sharing the same $LU$ factorization, i.e $(2+\frac{1}{3})$ products, to the total $\pi_{m}$.

\begin{figure}[!ht] 
\centering
  \includegraphics[width=\textwidth]{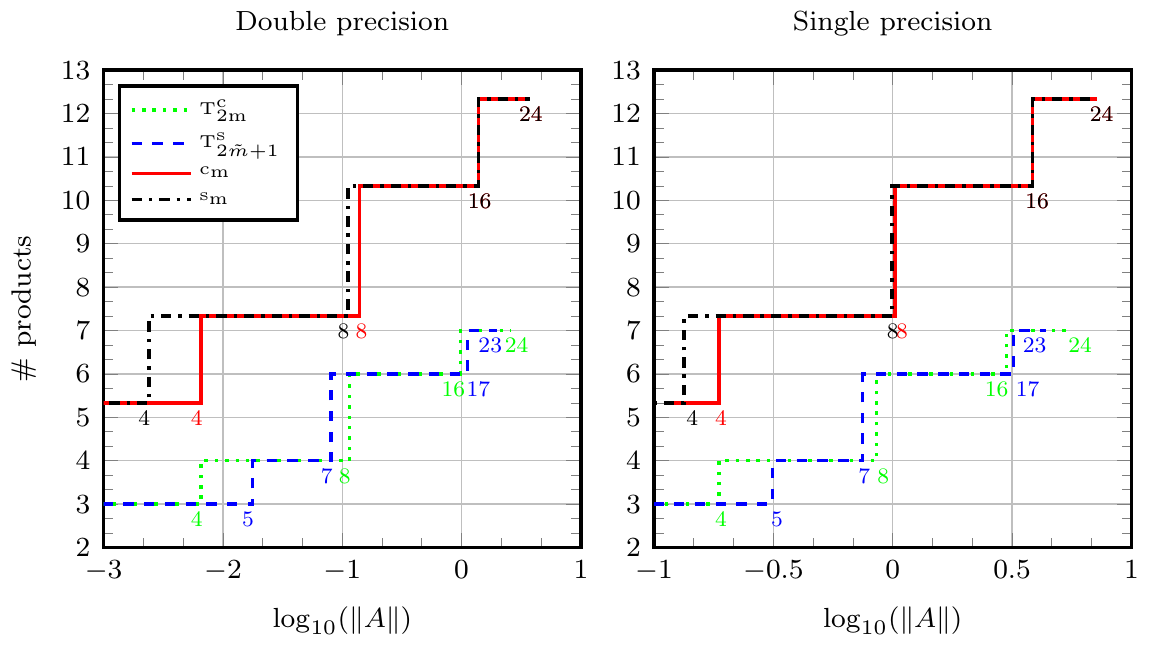}
\caption{\label{fig_stepfg} \small{Orders and corresponding number of products of each method versus $\left\|A\right\|$ in double and single precision floating-point arithmetic.} }
\end{figure}

The comparison of the theoretical performance of the new Taylor polynomial approximations $T_{2m}^{c},T_{2\tilde{m}+1}^{s}$  (with orders $ \left\{4, 8, 16, 24\right\}$ and $\left\{5,7,17,23\right\}$ respectively) and the Pad\'e approximations $c_m, s_m$ \cite{almohy14naf} (with orders $ \left\{4, 8, 16, 24\right\}$) has been illustrated in Figure \ref{fig_stepfg}: here
we plot $\left\|A\right\|$ versus the number of matrix products required for each approximation of  $\cos(A)$ and $\sin(A)$ simultaneously, both in
double (left) and single (right) precision. From the figure the improvement achieved by the proposed Taylor polynomial approximations 
is apparent.

\begin{table}[!ht]  \label{table1}
\caption{\small{Number of matrix multiplications $\Pi_{2m}$ and forward absolute error bounds $\theta_{2m}$ in double precision floating-point arithmetic, $u\leq 2^{-53}$, for the new Taylor algorithms $T_{2m}^{c}$, $T_{2\tilde{m}+1}^{s}$ and Pad\'e approximations $c_m$, $s_m$ \cite{almohy14naf}. The cost of the computation of two inverse matrices sharing the same $LU$ factorization, i.e $(2+\frac{1}{3})$, has been included in the cost $\pi_{m}$ for the Pad\'e approximations.} }\centering

\

\begin{tabular}{llllllll}
\hline
$2\tilde{m}$ & 4 & 6 & 16 & 22 \\ 
$2m$  & 4 & 8 & 16 & 24 \\ 
$\theta_{2m}^{c_m}$&6.5633e-3& 1.3959e-1& 1.3879 & 3.7288 \\
$\theta_{2m}^{s_m}$	&2.4019e-3&1.1213e-1& 1.3784 & 3.7287 \\
$\pi_{m}$ \cite{almohy14naf}	&  ${\bf 5+\frac{1}{3}}$ & ${\bf 7+\frac{1}{3}}$  & 
${\bf 10+\frac{1}{3}}$ & ${\bf 12+\frac{1}{3}}$ \\
$\theta_{2m}^{T_{2m}^{c}}$& 6.5633e-3&1.1495e-1 & 9.8108e-1&2.5675 \\ 
$\theta_{2m}^{T_{2\tilde{m}+1}^{s}}$&1.777e-2&8.0438e-2&1.1184&1.97  \\
$\Pi_{2m}$	&${\bf 3}$ & ${\bf 4}$& ${\bf 6}$& ${\bf 7}$ \\
\hline
\end{tabular} \label{doubleprctable}
\end{table}

\begin{table}[!ht]
\caption{\small{Same as Table \ref{table1}, but now in single precision floating-point arithmetic.} }\centering 

\

\begin{tabular}{llllllll}
 \hline 
$2\tilde{m}$ & 4 & 6 & 16 & 22 \\ 
$ 2m$  & 4 & 8 & 16 & 24 \\                  
$\theta_{2m}^{c_m}$	&1.8687e-1&1.0218& 3.8571 & 7.1575 \\
$\theta_{2m}^{s_m}$	&1.3355e-1&9.9511e-1& 3.8569 & 7.1575 \\  
$\pi_{m}$ \cite{almohy14naf}	&  ${\bf 5+\frac{1}{3}}$ & ${\bf 7+\frac{1}{3}}$  & 
${\bf 10+\frac{1}{3}}$ & ${\bf 12+\frac{1}{3}}$ \\
$\theta_{2m}^{T_{2m}^{c}}$&1.8709e-1 &8.5756e-1 &2.9935&5.5555\\ 
$\theta_{2m}^{T_{2\tilde{m}+1}^{s}}$&3.1386e-1 & 7.492e-1&3.2152 &4.3819  \\
$\Pi_{2m}$	&${\bf 3}$ & ${\bf 4}$& ${\bf 6}$& ${\bf 7}$ \\
\hline
\end{tabular} \label{singleprctable}
\end{table}

\begin{figure}[!ht] 
\centering
  \includegraphics[width=\textwidth]{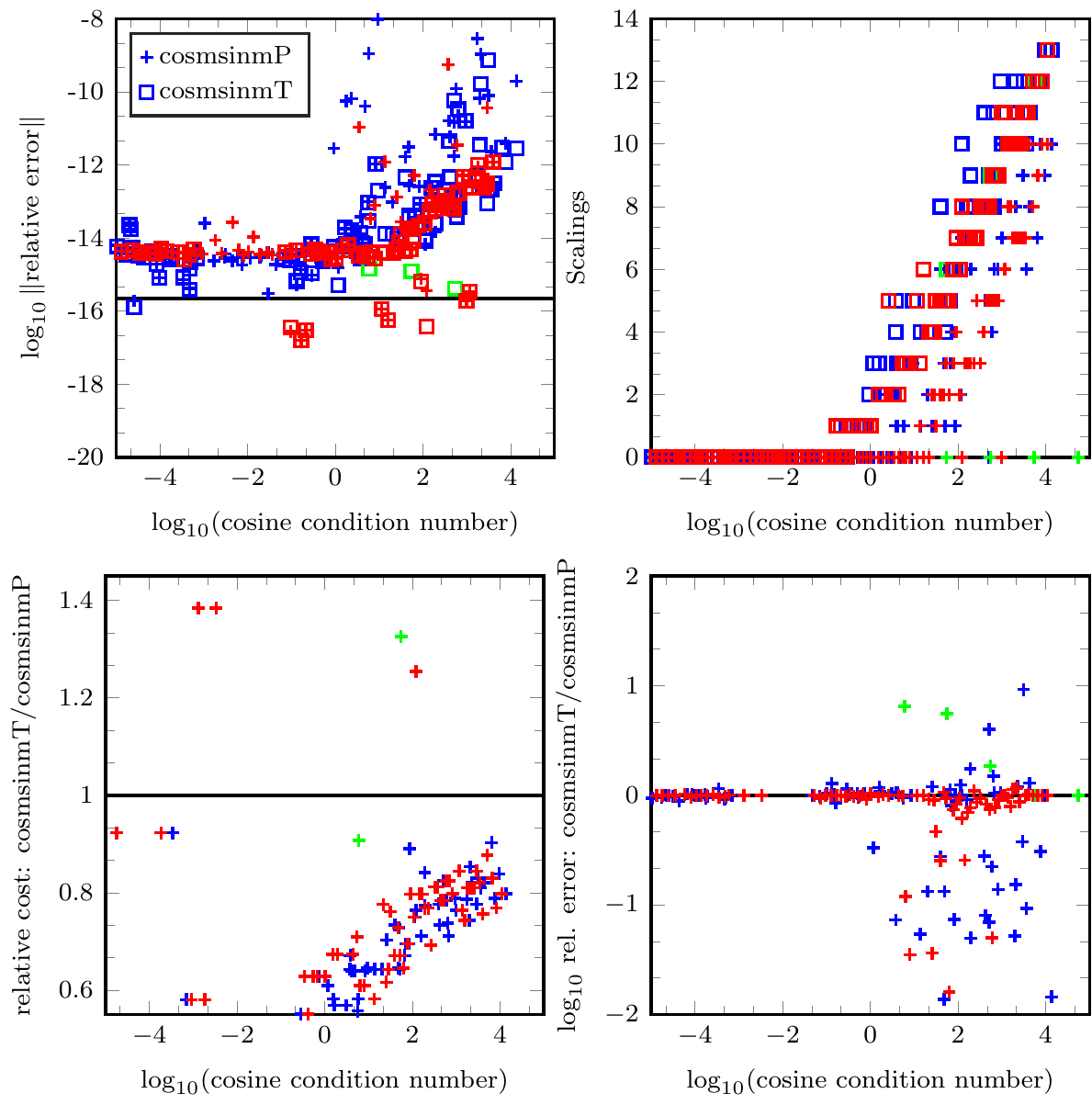}
\caption{  \small{Comparison of the results for the matrices of dimension $\leq 16\times16$ for cosine function}. \label{coserr_fg}}
\end{figure}

\begin{figure}[!ht] 
\centering
  \includegraphics[width=\textwidth]{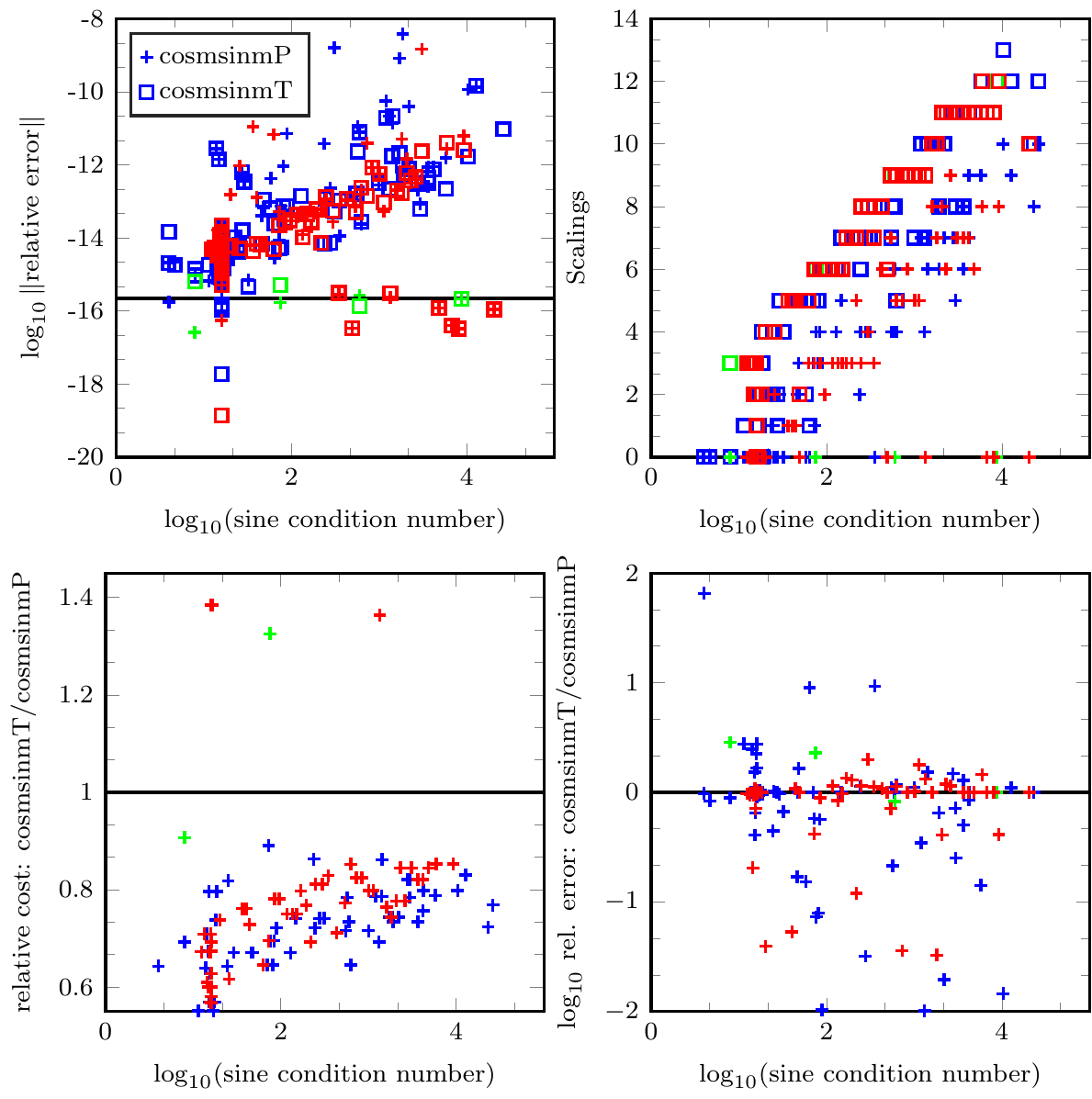}
\caption{\small{Comparison of the results for the matrices of dimension $\leq 16\times16$ for sine function.} \label{sinerr_fg}}
\end{figure}

\begin{figure}[!ht] 
\centering
  \includegraphics[width=\textwidth]{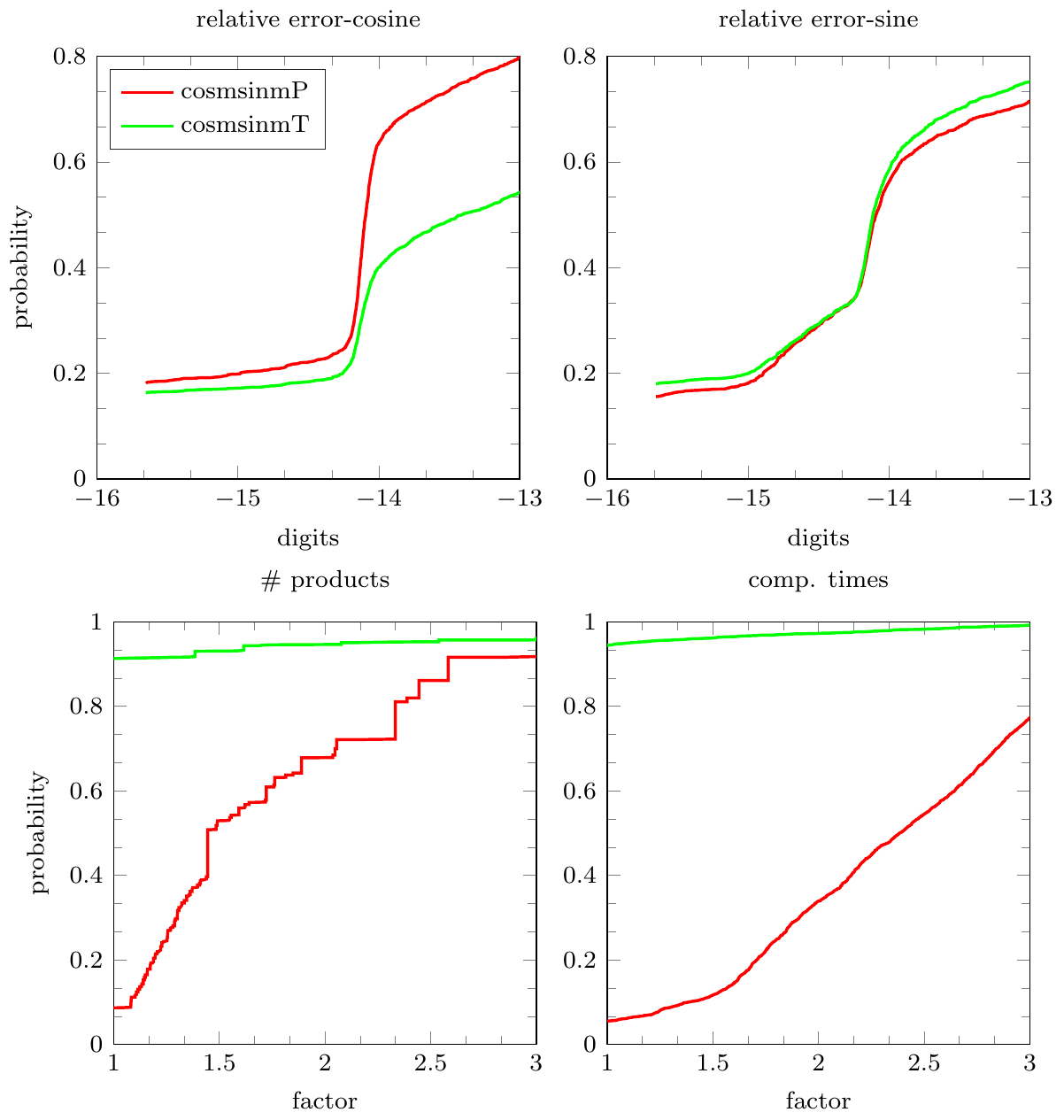}
\caption{\small{Performance profiles for the matrices of dimension $\leq 64\times64$. } \label{perf64_fg} }
\end{figure}

\begin{figure}[!ht] 
\centering
  \includegraphics[width=\textwidth]{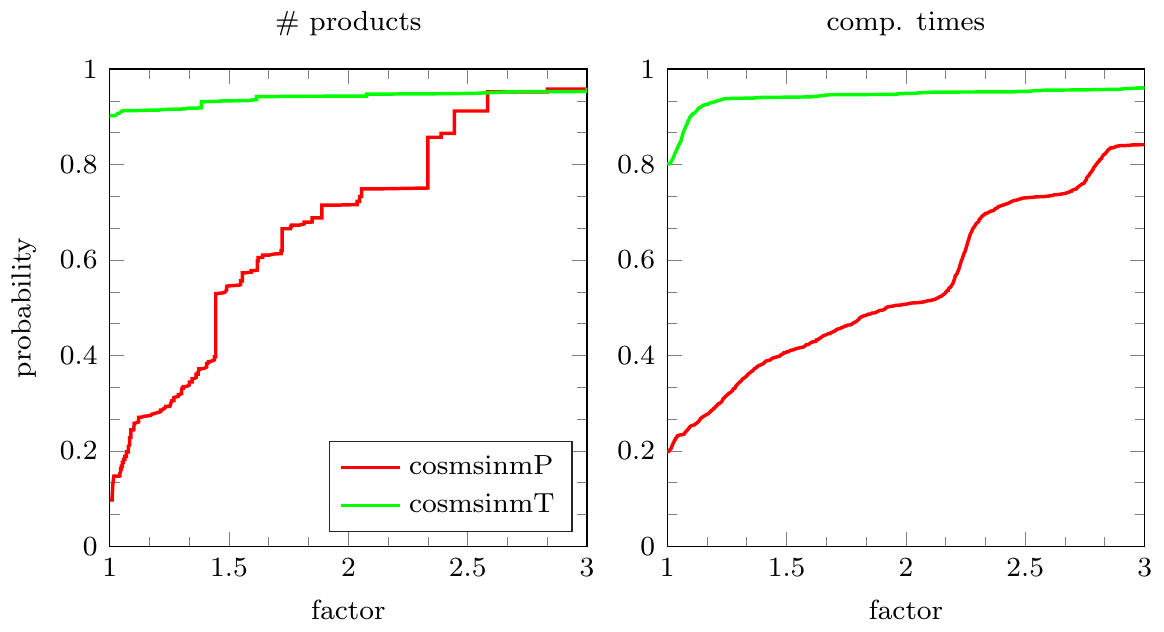}
\caption{\small{Performance profiles for the matrices of dimension $\leq 1024\times1024$}. \label{perf1024_fg} }
\end{figure}

\section{Numerical experiments}   \label{num.exp}

We measure the performance of new Taylor polynomials (denoted as `cosmsinmT') and the Pad\'e approximations (denoted as `cosmsinmP') \cite{almohy14naf} to compute matrix cosine and sine functions simultaneously. The platform of all numerical experiments is MATLAB R2013a and the matrix 1-norm has been used in implementing the algorithms. The experiments have been carried out for 2500 matrices (adjusted in order to have different norms) of the following cases:
\begin{itemize}
 \item 52 test matrices have been chosen from the MATLAB gallery function \cite{higham08fom} (blue). $690$ sampled matrices with different norms were tested.
\item Using \texttt{rand()} and \texttt{randn()} functions in MATLAB to randomly generate matrices with entries drawn from different distributions. $400$ matrices normally distributed, $500$ matrices uniformly distributed in the interval $(0,1)$ and $501$ matrices in the interval $(-0.5,0.5)$.     
\item Using \texttt{spdiags()} and \texttt{rand()} functions in MATLAB to construct 400 triangular nilpotent matrices with random rank (red).
\item $9$ matrices of the form
\begin{equation}\label{eq:symplectic2}
	A =
		\left( \begin{array}{cc} 
	1 & \lambda \\ 0 & -1
	\end{array}	\right),
\end{equation}
where $\lambda=1,10,\dots, 10^{8}$ (green), possibly leading to \textit{overscaling} (utilization of large value of scaling parameter $s$). 
\end{itemize}   
The same test matrices have been generated as in Remark 5 of \cite{bader19ctm} and all matrices are adjusted to have 1-norms over $(10^{-4}, 10^{4.1})$ in all numerical experiments. The condition numbers of each matrix function are computed by executing the function \texttt{funm$_{-}$condest1} from the Matrix Function Toolbox \cite{higham08fom}. The reference solutions of the matrix cosine and sine have been calculated with \textit{Mathematica} with 100 digits of precision. We have computed the relative error 
\[
    \frac{\left\|F-f(A)\right\|_{2}}{\left\|f(A)\right\|_{2}},
\]
where $F$ is an approximated value of $f(A)$.
In the following we show the results for double precision (similar results are obtained for single precision).
 We have simulated the results for the $2500$ matrices of dimension $\leq 16\times16$ in Figs.~\ref{coserr_fg}, \ref{sinerr_fg}. From the top left of the Figs.~\ref{coserr_fg}, \ref{sinerr_fg}, in general, the relative errors of both cosmsinmP and cosmsinmT methods produced in approximating the matrix cosine and sine functions change between $1.0e-12$ and $1.0e-15$ and they drop below the machine accuracy for few matrices. It can be observed from the top right of the Figs.~\ref{coserr_fg}, \ref{sinerr_fg}, the cosmsinmT method involves more scalings, particularly the cosmsinmP and cosmsinmT methods have leaded to the scaling for $741$ and $1180$ matrices respectively. Regarding to the bottom left of the Figs.~\ref{coserr_fg}, \ref{sinerr_fg}, the ratios of the cost $\text{cosmsinmT}/\text{cosmsinmP}$ are in general below 1, it also has been concluded from Tables \ref{doubleprctable}, \ref{singleprctable} and Fig.\ref{fig_stepfg}, the new method cosmsinmT requires less number of matrix products. As can be seen from the bottom right of the Figs.~\ref{coserr_fg}, \ref{sinerr_fg}, the accuracy of both methods is in good agreement with the theoretical results we have obtained. In these cases, some of the values of the relative errors have been replaced by machine accuracy (if these are lower) in the results of both methods. Furthermore, we plot performance profiles of the algorithms on a set of the test matrices exemplified for Figs.~\ref{coserr_fg}, \ref{sinerr_fg} in terms of the relative errors, number of products and computational times. The performance plot shows the percentage of problems (y-axis) that are within a given factor (x-axis) of the best method \cite{dolan02more}. In the experiments
illustrated by the performance profiles in Fig.~\ref{perf64_fg}, the $2500$ matrices of dimension $\leq 64\times64$ have been tested. We have observed that the cosmsinmT method has a lower relative error for 877 and 1257 of the 2500 matrices than the cosmsinmP method for computing the approximate values of the matrix cosine and sine functions respectively (358 and 348 results are equal). These results
are evident from the Fig.~\ref{perf64_fg} on the top.  It is seen clearly from the bottom of Fig.~\ref{perf64_fg} that the cosmsinmT method is less expensive than cosmsinmP.

The performance profiles in Fig. \ref{perf1024_fg} resulted from demonstrating the returns from the $2500$ matrices of dimension $\leq 1024\times1024$ confirm the superiority of the cosmsinmT method in the sense of computational cost.

\section{Conclusions}

We have presented a new algorithm to compute the matrix cosine and sine. The algorithm contains several methods that are optimised for different values of the norm of the matrix and the desired accuracy, and can be combined with the scaling  and squaring technique. Each of these methods is obtained by following a sequence in which each stage uses the results from all previous ones. An error analysis is also carried out and we have shown both theoretically as well as in the numerical experiments that the new algorithm is superior to other procedures from the literature that are based on Pad\'e approximations to the  matrix cosine and sine. 

The new algorithm only involves matrix-matrix products and does not require to compute the inverse of matrices as it the the case of the Pad\'e approximations. The cost to compute the inverse of a dense matrix can be taken as 4/3 the cost of the product of two dense matrices. However, for sparse matrices, the computational cost of the proposed algorithms grow nearly linearly while the cost of Pad\'e approximations grows much faster because, in general, the inverse of a sparse matrix is a dense matrix. 

\subsection*{Acknowledgments}
The work of MS has been funded by The Scientific and Technological Research Council of Turkey (TUBITAK) with Grand Number 1059B191802292. 
PB, SB and FC acknowledge financial support from Ministerio de Econom\'{\i}a, Industria y Competitividad (Spain) through projects MTM2016-77660-P 
and PID2019-104927GB-C21 (AEI/FE\-DER, UE).


\begin{thebibliography}{00}

\bibitem{almohy14naf}
 A.H. Al-Mohy, N.J. Higham and S.D. Relton,
New algorithms for computing the matrix sine and cosine separately or simultaneously, 
SIAM J. Sci. Compu.  37 (2015)  A456 - A487.

\bibitem{preprint} 
{P. Bader, S. Blanes, and F. Casas,}
 An improved algorithm to compute the exponential of a matrix, arXiv:1710.10989 [math.NA]  (2017) preprint.

 
\bibitem{bader19ctm} 
{P. Bader, S. Blanes, and F. Casas,}
Computing the matrix exponential with an optimized Taylor polynomial approximation,
Mathematics 7 (2019) 1174 doi:10.3390/math7121174.

\bibitem{bader17sif} 
P. Bader, S. Blanes, E. Ponsoda, and M Seydao\u{g}lu,
Symplectic integrators for the matrix Hill's equation and its applications to engineering models,
J. Comput. Appl. Math. 316 (2017) 47 - 59.

\bibitem{bader15tss} 
{P. Bader, S. Blanes, and M. Seydao\u{g}lu,}
 The scaling, splitting and squaring method for the exponential of perturbed matrices,
SIAM  J. Matrix Anal. Appl.  36 (2015) 594 - 614.

 \bibitem{blanes16aci} 
S. Blanes, and F. Casas, 
A Concise Introduction to Geometric Numerical Integration,
CRC Press: Boca Raton, FL, USA, 2016.

\bibitem{blanes15aea} 
S. Blanes,  F. Casas, and A. Murua, 
An efficient algorithm based on splitting for the time integration of the Schr\"odinger equation,
 J. Comput. Phys.  303 (2015) 396 - 412.

\bibitem{blanes09tme}
S. Blanes, F. Casas, J.A. Oteo, and J. Ros,
The {M}agnus expansion and some of its applications,
Phys. Rep.  470 (2009) 151 - 238.

\bibitem{blanes02hoo}  
S. Blanes, F. Casas and J. Ros, 
High order optimized geometric integrators for linear differential equations, BIT, 42 (2002) 262 - 284.

\bibitem{hairer06gni} 
E. Hairer, C. Lubich, and G. Wanner, 
 Geometric Numerical Integration. Structure-Preserving Algorithms for Ordinary Differential Equations, 2nd Ed., Springer, Berlin, 2006.


\bibitem{higham08fom} 
N.J. Higham, 
Functions of Matrices: Theory and Computation,
Society for Industrial and Applied Mathematics, Philadelphia, PA, USA, 2008.

\bibitem{higham10cma} 
N.J. Higham, and A.H. Al-Mohy,
Computing matrix functions,
 Acta Numerica   19 (2010) 159 - 208.

\bibitem{hochbruck97oks}
M. Hochbruck and C. Lubich, 
On Krylov subspace approximations to the matrix exponential operator, 
SIAM J. Numer. Anal.  34 (1997) 1911 – 1925. 

\bibitem{hochbruck10ein} 
M. Hochbruck and  A. Ostermann,
Exponential integrators, 
Acta Numerica  19 (2010) 209 - 286.


\bibitem{iserles00lgm1} 
A. Iserles,
A First Course in the Numerical Analysis of Differential Equations,
Cambridge University Press, 2nd ed.,
2008

\bibitem{iserles00lgm} 
A. Iserles, H.Z. Munthe-Kaas, S.P. N{\o}rsett, and A. Zanna,
Lie group methods, 
Acta Numerica  9 (2000)215 - 365.

\bibitem{lei92afa} L. Lei, and T. Nakamura,
A fast algorithm  for  evaluating the matrix polynomial $I+A+ \cdots + A^{N-1}$,
IEEE Trans. Circuits Sys.-I:  Fund. Theory Appl. 39 (1992) 299 - 300.

\bibitem{leimkuhler04shd} 
B. Leimkuhler and S. Reich, 
Simulating Hamiltonian Dynamics, 
Cambridge University Press, 2004.


\bibitem{liang04fmf} 
W. Liang, R. Baer, C. Saravanan, Y. Shao, A.T. Bell, M. Head-Gordon,
Fast methods for resumming matrix polynomials and Chebyshev matrix polynomials,
J. Comput. Phys. 194 (2004) 575 - 587.


\bibitem{lubich08fqt} 
C. Lubich, 
From Quantum to Classical Molecular Dynamics: Reduced Models and Numerical Analysis,
European Mathematical Society, 2008.

\bibitem{moler03ndw} 
C.B. Moler, and C.F. Van Loan,
Nineteen dubious ways to compute the exponential of a matrix, twenty-five years later,
SIAM Review 45 (2003) 3 - 49.


\bibitem{paterson73otn} 
M.S. Paterson, and L.J. Stockmeyer,
On the number of nonscalar multiplications necessary to evaluate polynomials,
SIAM J. Comput.  2 (1973) 60 - 66.

\bibitem{sanzserna94nhp} 
{J.M. Sanz-Serna and M.P. Calvo}, 
 Numerical Hamiltonian Problems,  
Chapman \& Hall, London, 1994.

\bibitem{sastre18eeo} 
J. Sastre,
Efficient evaluation of matrix polynomials,
 Linear Algebra Appl.539 (2018) 229 - 250.

\bibitem{sastre19ftp}
J. Sastre, J. Ib\'a\~nez, P. Alonso-Jord\'a, J. Peinado, E. Defez,
Fast Taylor polynomial evaluation for the computation of the matrix cosine,
J. Comput. Appl. Math. 354 (2019) 641 - 650.

\bibitem{sastre19btc}
J. Sastre, J. Ibáñez, and E. Defez,
Boosting the computation of the matrix exponential,
Appl. Math. Comput. 340 (1019) 206 - 220.

\bibitem{sidje98eas} 
R.B. Sidje, 
Expokit: a software package for computing matrix exponentials,
ACM Trans. Math. Software   24 (1998) 130 - 156.


\bibitem{vanloan79ano} Van Loan, C.
A note on the evaluation of matrix polynomials,
IEEE Transactions on Automatic Control 24 (1979) 320 - 321.


\bibitem{westreich89etm} Westreich, D.
Evaluating  the  matrix polynomial $I+A+ \cdots + A^{N-1}$,
IEEE Trans. Circuits Sys.  36 (1989) 162 - 164.

\bibitem{dolan02more} E. D. Dolan and J. J. More, Benchmarking optimization software with performance profiles,
Math. Programming 91 (2002) 201 - 213.



\end{thebibliography}
\end{document}